\documentclass[11pt,letterpaper,reqno]{amsart}
\usepackage{amsmath,amsthm,amsfonts,amssymb,mathtools,algpseudocode,epsfig}
\usepackage{comment,bookmark,bm,color}
\usepackage{caption,subcaption}

\hypersetup{pdfstartview={FitH}}

\newtheorem{Theorem}{{\bf Theorem}}[section]

\newtheorem{Corollary}[Theorem]{{\bf Corollary}}
\newtheorem{Algorithm}{{\bf Algorithm}}

\numberwithin{equation}{section}

\newcommand{\calX}{\mathcal{X}}
\newcommand{\calY}{\mathcal{Y}}
\newcommand{\calA}{\mathcal{A}}
\newcommand{\calB}{\mathcal{B}}
\newcommand{\calS}{\mathcal{S}}
\newcommand{\calG}{\mathcal{G}}
\newcommand{\calE}{\mathcal{E}}
\newcommand{\C}{\mathbb{C}}

\begin{document}

\title[Hybrid CUR-type decomposition]{Hybrid CUR-type decomposition of tensors in the Tucker format}

\author{Erna Begovi\'{c}~Kova\v{c}}\thanks{Erna Begovi\'{c} Kova\v{c}, Faculty of Chemical Engineering and Technology, University of Zagreb, Maruli\'{c}ev trg 19, 10000 Zagreb, Croatia, \texttt{ebegovic@fkit.hr}}
\date{\today}

\renewcommand{\subjclassname}{\textup{2020} Mathematics Subject Classification}

\begin{abstract}
The paper introduces a hybrid approach to the CUR-type decomposition of tensors in the Tucker format.
The idea of the hybrid algorithm is to write a tensor $\calX$ as a product of a core tensor $\calS$, a matrix $C$ obtained by extracting mode-$k$ fibers of $\calX$, and matrices $U_j$, $j=1,\ldots,k-1,k+1,\ldots,d$, chosen to minimize the approximation error. The approximation can easily be modified to preserve the fibers in more than one mode.
The approximation error obtained this way is smaller than the one from the standard tensor CUR-type method.
This difference increases as the tensor dimension increases. It also increases as the number of modes in which the original fibers are preserved decreases.
\end{abstract}

\maketitle

\section{Introduction}

Let $A\in\mathbb{C}^{m\times n}$ be a matrix with singular value decomposition given by
$$A=U\Sigma V^*=\left[
                 \begin{array}{cc}
                   U_{k} & U_{m-k} \\
                 \end{array}
               \right]\left[
                 \begin{array}{cc}
                   \Sigma_k & 0 \\
                   0 & \Sigma_{\min\{m,n\}-k} \\
                 \end{array}
               \right]\left[
                 \begin{array}{cc}
                   V_{k} & V_{n-k} \\
                 \end{array}
               \right]^*.
$$
Its best rank-$k$, $k<\min\{m,n\}$, approximation is obtained using truncated SVD,
\begin{equation}\label{svdapproximation}
A_k=U_k\Sigma_k V_k^*.
\end{equation}
However, approximation~\eqref{svdapproximation} is often hard to interpret in applications, especially when we work with very big matrices. Besides, this approximation does not keep useful matrix properties like sparsity and non-negativity. Therefore, in recent years attention is given to low-rank approximations obtained by interpolatory factorizations. These approximations are suboptimal, but keep the properties mentioned above and are more suitable in some applications where the columns and/or rows should keep their original meaning.

The best known examples of the interpolatory factorizations are CX and CUR factorizations.
A matrix CX factorization of $A\in\C^{m\times n}$ is a decomposition of the form
\begin{equation}\label{cx}
A=CX,
\end{equation}
where $C\in\C^{m\times k}$ contains $k$ columns of $A$ and $X\in\C^{k\times n}$.
A matrix CUR factorization of $A\in\C^{m\times n}$ is decomposition of the form
\begin{equation}\label{cur}
A=CUR,
\end{equation}
where $C\in\C^{m\times k}$ contains $k$ columns of $A$, $R\in\C^{k\times n}$ contains $k$ rows of $A$ and $X\in\C^{k\times k}$.
In these factorizations columns of $C$ and rows of $R$ keep the original interpretation from $A$.
Usually, we are not looking for the exact decompositions~\eqref{cx} or~\eqref{cur}, but for the approximations
$$A=CX+E \quad \text{or} \quad A=CUR+E, \quad \text{where } \|E\|\ll\|A\|.$$

Interpolatory decomposition of a tensor $\calX\in\C^{n_1\times n_2\times\cdots\times n_d}$ in the Tucker representation is the product of a core tensor $\calG\in\C^{r_1\times r_2\times\cdots\times r_d}$ and matrices $C_j\in\C^{n_j\times r_j}$, $1\leq j\leq d$.
Over the last ten years different authors have studied generalizations of the CUR decomposition. In~\cite{MMD08} the authors study CUR tensor decomposition where the fibers that define the decomposition are chosen randomly according to a specific probability distribution.
In~\cite{CC10} an adaptive algorithm that sequentially selects fibers of a $3$rd order tensor $\calX$ that form matrices $C_n$, $n=1,2,3$, is developed.
Different choices for matrices $C_n$ are studied in~\cite{Saibaba16}, where the author gives detailed analysis of the computational costs and error bounds.

\textbf{Motivation for the hybrid approach:}
Generalizations of CUR decomposition for tensors represented in the Tucker format give the approximation error that depends on the dimension $d$. This can be a problem for high-dimensional tensors. On the other hand, in the applications it is not always important to keep the original entries in all modes. In the same way as the matrix $CX$ decomposition preserves only the original columns of a starting matrix, in the tensor case we can keep the original fibers in only one mode, or in more, but not all modes.
The idea of the hybrid algorithm is to write a tensor $\calX\in\C^{n_1\times n_2\times\cdots\times n_d}$ as a product of a core tensor $\calS\in\C^{r_1\times r_2\times\cdots\times r_d}$, a matrix $C\in\C^{n_k\times r_k}$ obtained by extracting mode-$k$ fibers of $\calX$, and matrices $U_j\in\C^{n_j\times r_j}$, $j=1,\ldots,k-1,k+1,\ldots,d$, chosen to minimize the approximation error.
This difference between the error obtained by the hybrid approach and the error from the tensor CUR method gets more important as the tensor dimension increases. We keep the Tucker representation because it is one of the most commonly used tensor formats. Also, it is suitable for function-related tensors.

In Section~\ref{sec:preliminaries} we introduce the notation and give an overview of the concepts used in the paper. We introduce the hybrid method in Section~\ref{sec:hybrid}. Moreover, we give the error bound for the new method and compare it to the error resulting from the standard CUR approach. In Section~\ref{sec:numerical} we present the results of several numerical tests. We refer to the matrix case in Section~\ref{sec:matrix}.

\section{Preliminaries and notation}\label{sec:preliminaries}

Throughout the paper we use tensor notation from~\cite{KB09}. The \emph{order} of a tensor is a number of its dimensions. Tensors  of order three or higher are denoted by calligraphic letters, e.g.\@ $\calX$, while matrices (order two tensors) are, as usually, denoted by capital letters, e.g.\@ $A$.
Tensor analogues of rows and columns are called \emph{fibers} and they are extracted from a tensor by fixing all indices but one. For a third order tensor its fibers are columns, rows, and tubes, denoted $x_{:jk},x_{i:k},x_{ij:}$, respectively.
If we fix all but two indices, we get tensor \emph{slices}. In the case of a third order tensor the slices are matrices $X_{i::},X_{:j:},X_{::k}$. The \emph{norm} of a tensor $\calX\in\C^{n_1\times n_2\times\cdots\times n_d}$ is a generalization of the matrix Frobenius norm and it is given by relation
$$\|\calX\|_F=\sqrt{\sum_{i_1=1}^{n_1}\sum_{i_2=1}^{n_2}\cdots\sum_{i_d=1}^{n_d} |x_{i_1i_2\ldots i_d}|^2}.$$

Tensor \emph{unfolding} is a reordering of an order-$d$ tensor into a matrix. The mode-$m$ unfolding, $1\leq m\leq d$, of $\calX\in\C^{n_1\times n_2\times\cdots\times n_d}$ is an $n_m\times(n_1\cdots n_{m-1}n_{m+1}\cdots n_d)$ matrix $X_{(m)}$ obtained by arranging mode-$m$ fibers of $\calX$ into columns of $X_{(m)}$.

The \emph{mode-$m$ product} of a tensor $\calX\in\C^{n_1\times n_2\times\cdots\times n_d}$ with a matrix $U\in\C^{p\times n_m}$ is a $n_1\times\cdots\times n_{m-1}\times p\times n_{m+1}\times\cdots\times n_d$ tensor
\begin{equation}\label{eq:product}
\calY=\calX\times_mU \quad \text{such that} \quad Y_{(m)}=UX_{(m)}.
\end{equation}
Elementwise, relation~\eqref{eq:product} can be written as
$$\calY_{i_1\ldots i_{m-1}ji_{m+1}\ldots i_d}=\sum_{i_m=1}^{n_d} x_{i_1i_2\ldots i_d}u_{ji_m}.$$
We will use the associativity properties of mode-$m$ product,
\begin{equation}\label{eq:associativity}
\begin{aligned}
\calX\times_m M \times_m N & =\calX\times_m(NM) \quad \text{and} \\
\calX\times_m M\times_n N & =\calX\times_n N\times_m M, \quad \text{for} \ m\neq n.
\end{aligned}
\end{equation}

\emph{Tucker decomposition} is a decomposition of a tensor $\calX$ into a core tensor $\calS$ multiplied by a matrix in each mode,
$$\calX=\calS\times_1U_1\times_2U_2\times_3\cdots\times_dU_d.$$
If $U_1,U_2,\ldots,U_d$ are unitary matrices, this decomposition is referred to as higher order singular value decomposition (HOSVD) studied in~\cite{DeL-hosvd}. Matrices $U_j$, $1\leq j\leq d$, from HOSVD are computed using the SVD of each unfolding of $\calX$.
Note that if
$$\calY=\calX\times_1M_1\times_2M_2\times_3\cdots\times_dM_d,$$
for some matrices $M_i$, $1\leq i\leq d$ of the adequate size, then
\begin{equation}\label{eq:productmatrix}
Y_{(m)}=M_mX_{(m)}\big{(}M_d\otimes M_{d-1}\otimes\cdots\otimes M_{m+1}\otimes M_{m-1}\otimes\cdots\otimes M_1\big{)}^*, \quad  \text{for } 1\leq m\leq d.
\end{equation}

\emph{Multilinear rank} of $\calX$ is a $d$-tuple $(r_1,r_2,\ldots,r_d)$, where $r_j=\text{rank}(X_{(j)})$, $1\leq j\leq d$.
If $r_1=\cdots=r_d=k$ we say that $\calX$ is a rank-$k$ tensor.
A simple way to obtain a rank-$(r_1,r_2,\ldots,r_d)$ approximation $\hat{\calX}$ of $\calX$ is using truncated HOSVD (T-HOSVD), which is a tensor analogue of~\eqref{svdapproximation}. It is described in Algorithm~\ref{alg:thosvd}.

\begin{Algorithm}\label{alg:thosvd}
\hrule\vspace{1ex}
\emph{T-HOSVD}
\vspace{0.5ex}\hrule\vspace{0.5ex}
\begin{algorithmic}
\For{$i=1,\ldots,d$}
\State Compute matrix $U_i$ containing the leading $r_i$ left singular vectors of $X_{(i)}$.
\EndFor
\State $\calS=\calX \times_1 U_1^* \cdots\times_d U_d^*$
\State $\hat{\calX}=\calS\times_1U_1 \cdots \times_d U_d$
\end{algorithmic}
\hrule
\end{Algorithm}

The core tensor from HOSVD is, in general, not diagonal. Thus, HOSVD does not lead to the best low multilinear rank tensor approximation, as it is the case with SVD.
To improve the approximation obtained by T-HOSVD one can use an iterative algorithm with the initialization based on the result obtained by T-HOSVD.
A popular iterative algorithm for low-rank tensor approximation is the higher-order orthogonal iteration (HOOI)~\cite{SaadHOOI2006}. If the starting tensor is symmetric or anti-symmetric, a good choice can be the structure-preserving Jacobi algorithm~\cite{Ishteva13,BegKre17}.

\subsection{Tensor CUR decomposition}

A CUR-type decomposition of a tensor $\calA\in\C^{n_1\times n_2\times\cdots\times n_d}$ is given by
\begin{equation}\label{cur-tensor}
\calA\approx\calS\times_1C_1\times_2C_2\times_3\cdots\times_dC_d,
\end{equation}
where $\calS\in\C^{r_1\times r_2\times\cdots\times r_d}$ is a core tensor and matrices $C_j\in\C^{n_j\times r_j}$, $1\leq j\leq d$, contain $r_j$  mode-$j$ fibers of $\calA$. The algorithm that we are going to use for the CUR-type tensor decomposition is higher order interpolatory decomposition (HOID) from~\cite{Saibaba16}. It is derived for the tensors in the Tucker format and it shows good numerical behavior. This decomposition is based on CX decompositions of the unfoldings $A_{(j)}$ of $\calA$.

One way to compute matrices $C_j$, $1\leq j\leq d$, from~\eqref{cur-tensor} is using QR decomposition with column pivoting (PQR) of $A_{(j)}$,
\begin{equation}\label{rel:PQR}
A_{(j)}P=QR, \quad 1\leq j\leq d,
\end{equation}
where $P$ is a permutation matrix, $Q$ is unitary and $R$ is upper-triangular.
We can write relation~\eqref{rel:PQR} using block partitions as
\begin{equation}\label{rel:QR}
A_{(j)}\left[
           \begin{array}{cc}
             P_1 & P_2 \\
           \end{array}
         \right]=\left[
           \begin{array}{cc}
             Q_1 & Q_2 \\
           \end{array}
         \right]\left[
                  \begin{array}{cc}
                    R_{11} & R_{12} \\
                    0 & R_{22} \\
                  \end{array}
                \right], \quad 1\leq j\leq d,
\end{equation}
where $P_1$ and $Q_1$ have $r_j$ columns and $R_{11}$ is $r_j\times r_j$.
Then we set
\begin{equation}\label{rel:C}
C_j=Q_1R_{11},
\end{equation}
and
\begin{equation}\label{rel:CX}
A_{(j)}=C_jF+E,
\end{equation}
where
\begin{equation}\label{rel:CXerror}
F=\left[
           \begin{array}{cc}
             I & R_{11}^{-1}R_{12} \\
           \end{array}
         \right]P^T, \quad E=\left[
           \begin{array}{cc}
             0 & Q_2R_{22} \\
           \end{array}
         \right]P^T.
\end{equation}
Knowing $C_j$, $1\leq j\leq d$, the core tensor $\calS$ from~\eqref{cur-tensor} is obtained as
$$\calS=\calA \times_1 C_1^{+} \cdots\times_d C_d^{+}.$$
Here $C_j^{+}$ stands for the Moore-Penrose inverse of $C_j$, $1\leq j\leq d$.

HOID Algorithm for a tensor $\calA\in\C^{I_1\times I_2\times\cdots\times I_d}$ is presented in~\ref{alg:hoid}.

\begin{Algorithm}\label{alg:hoid}
\hrule\vspace{1ex}
\emph{HOID}
\vspace{0.5ex}\hrule\vspace{0.5ex}
\begin{algorithmic}
\For{$i=1,\ldots,d$}
\State Compute $C_i$ using PQR of $A_{(i)}$.
\EndFor
\State $\calS=\calA \times_1 C_1^{+} \cdots\times_d C_d^{+}$
\State $\hat{\calA}=\calS\times_1C_1 \cdots \times_d C_d$
\end{algorithmic}
\hrule
\end{Algorithm}

\section{Hybrid algorithm}\label{sec:hybrid}

In this section we describe and analyze our hybrid approach to tensor CUR-type decomposition.

Let $\calA\in\mathbb{C}^{n_1\times n_2\times\cdots\times n_d}$. We are looking for a low multilinear rank approximation $\hat{\calA}$ of $\calA$.  Precisely, we are looking for multilinear rank $(r_1,r_2,\ldots,r_d)$ tensor
\begin{equation}\label{def:hybrid-approx}
\hat{\calA}=\calS\times_1U_1\times_2\cdots\times_{m-1}U_{m-1}\times_mC_m\times_{m+1}U_{m+1}\cdots\times_dU_d,
\end{equation}
such that
$$\calA=\hat{\calA}+\calE, \quad \|\calE\|\ll\|\calA\|.$$
In~\eqref{def:hybrid-approx}, $n_m\times r_m$ matrix $C_m$ contains $r_m$ columns of $A_{(m)}$, that is mode-$m$ fibers of $\calA$, $\calS$ is a $r_1\times r_2\times\cdots\times r_d$ tensor, and $U_i$ are $n_i\times r_i$ matrices, $i=1,\ldots,m-1,m+1,\ldots,d$.
Without loss of generality we can assume that $m=1$. Thus, relation~\eqref{def:hybrid-approx} reads
\begin{equation}\label{def:hybrid-approx1}
\hat{\calA}=\calS\times_1C\times_2U_2\times_3\cdots\times_dU_d.
\end{equation}
We determine matrix $C$ using PQR decomposition as in relation~\eqref{rel:C} and as in Algorithm~\ref{alg:hoid}. On the other hand, to find $U_i$ we use SVD of $A_{(i)}$, $2\leq i\leq d$, as in Algorithm~\ref{alg:thosvd}. Then, the core tensor $\calS$ is obtained as
\begin{equation}\label{def:hybrid-core}
\calS=\calA\times_1C^{+}\times_2U_2^*\times_3\cdots\times_dU_d^*.
\end{equation}

Let us check that equation~\eqref{def:hybrid-core} gives an optimal $\calS$.
We are looking for the core tensor $\calS$ from~\eqref{def:hybrid-approx1} such that
\begin{equation}\label{eq:min}
\|\calA-\calS\times_1C\times_2U_2\times_3\cdots\times_dU_d\|_F\rightarrow\min.
\end{equation}
Using mode-$1$ matricizations and equation~\eqref{eq:productmatrix} we can write minimization problem~\eqref{eq:min} as
$$\|A_{(1)}-CS_{(1)}\big{(}U_d\otimes\cdots\otimes U_2\big{)}^*\|\rightarrow\min.$$
With the assumption that $C$ has full column rank, it follows from~\cite{FT07} that the optimal $S_{(1)}$ is
$$S_{(1)}=C^+A_{(1)}\Big{(}\big{(}U_d\otimes\cdots\otimes U_2\big{)}^*\Big{)}^*=C^+A_{(1)}\big{(}U_d\otimes\cdots\otimes U_2\big{)}.$$
Now we use~\eqref{eq:productmatrix} to go back on the tensor format. Thus, we get $\calS$ as in relation~\eqref{def:hybrid-core}.

The idea of the hybrid approach is summarized in Algorithm~\ref{alg:hybrid}. This algorithm corresponds to the problem given in~\eqref{def:hybrid-approx1}. It can easily be modified to $m=2,3,\ldots,d$.

\begin{Algorithm}\label{alg:hybrid}
\hrule\vspace{1ex}
\emph{Hybrid algorithm}
\vspace{0.5ex}\hrule\vspace{0.5ex}
\begin{algorithmic}
\For{$i=2,\ldots,d$}
\State Compute matrix $U_i$ containing the leading $r_i$ left singular vectors of $A_{(i)}$.
\EndFor
\State Compute $C$ using PQR of $A_{(1)}$.
\State $\calS=\calA\times_1C^{+}\times_2U_2^*\cdots\times_dU_d^*$
\State $\hat{\calA}=\calS\times_1C\times_2U_2\cdots\times_dU_d$
\end{algorithmic}
\hrule
\end{Algorithm}

In Algorithm~\ref{alg:hybrid} we can choose to extract fibers from more than one mode of $\calA$. Assume that the approximation of $\calA$ requires keeping fibers from the first $t$ modes. In this case we apply HOSVD to find $U_i$, $t+1\leq i\leq d$, and PQR decomposition to find $C_j$, $1\leq j\leq t$. Then we set
\begin{equation}\label{def:hybrid-more}
\begin{aligned}
\calS & =\calA\times_1C_1^{+}\cdots\times_tC_t^{+}\times_{t+1}U_{t+1}^*\cdots\times_dU_d^*, \\
\hat{\calA} & =\calS\times_1C_1\cdots\times_tC_t\times_{t+1}U_{t+1}\cdots\times_dU_d.
\end{aligned}
\end{equation}

Note that instead of using HOID Algorithm to compute matrix $C$, one can choose a different method, such as volume maximization in tensor approximations~\cite{OST08}, leverage scores method~\cite{DMM08,MMD08,MD09,BW17}, and discrete empirical interpolation method (DEIM)~\cite{DG16,SE16}. The latter two are extended to the tensor case in~\cite{Saibaba16}.

\subsection{Error analysis}

In Theorem~\ref{tm:error} we give the error bound for the low multilinear rank approximation~\eqref{def:hybrid-approx1}.

\begin{Theorem}\label{tm:error}
Let $\calA\in\mathbb{C}^{n_1\times n_2\times\cdots\times n_d}$. Let $\hat{\calA}$ be an approximation of $\calA$ computed by Algorithm~\ref{alg:hybrid}. Then the approximation error $\calE$ satisfies the following inequality,
\begin{equation}\label{rel:tmerror}
\|\calE\|_F^2=\|\calA-\hat{\calA}\|_F^2 \leq p(r_1,n_1)(n_1-r_1)\sigma_{r_1+1}^2(A_{(1)}) + \sum_{j=2}^d (n_j-r_j)\sigma_{r_j+1}^2(A_{(j)}),
\end{equation}
where
\begin{equation}\label{def:p}
p(r,n):=\left(1+2r+\sum_{j=1}^{r-1}4^j(r-j)\right)(n-r),
\end{equation}
and $\sigma_i(X)$ stands for the $i$th singular value of $X$.
\end{Theorem}

\begin{proof}
From~\eqref{def:hybrid-approx1} and~\eqref{def:hybrid-core}, using the properties of the mode-$m$ product given in~\eqref{eq:associativity}, we have
\begin{align*}
\|\calE\|_F^2 & =\|\calA-\calS\times_1C\times_2U_2\times_3\cdots\times_dU_d\|_F^2 \\
& = \|\calA-(\calA\times_1C^{+}\times_2U_2^*\times_3\cdots\times_dU_d^*)\times_1C\times_2U_2\times_3\cdots\times_dU_d\|_F^2 \\
& = \|\calA-\calA\times_1(CC^{+})\times_2(U_2U_2^*)\times_3\cdots\times_d(U_dU_d^*)\|_F^2.
\end{align*}
Matrices $CC^{+}$ and $U_jU_j^*$, $2\leq j\leq d$, are projections. Therefore, we can use the following result from~\cite[Lemma 2.1]{Saibaba16}:
$$\|\calX-(\calX\times_1\Pi_1\times_2\Pi_2\times_3\cdots\times_d\Pi_d)\|_F^2\leq\sum_{j=1}^d\|\calX-\calX\times_j\Pi_j\|_F^2,$$
that holds for projections $\Pi_1,\ldots,\Pi_d$.
It follows that
\begin{align}
\|\calE\|_F^2 & \leq \|\calA-\calA\times_1(CC^{+})\|_F^2+\sum_{j=2}^d\|\calA-\calA\times_j(U_jU_j^*)\|_F^2 \nonumber \\
& = \|(I_{n_1}-CC^+)A_{(1)}\|_F^2+\sum_{j=2}^d\|(I_{n_j}-U_jU_j^*)A_{(j)}\|_F^2. \label{tm:error1}
\end{align}

For $r<n$, set
$$\tilde{I}_{r,n}:=
\left[
  \begin{array}{cc}
    I_{r} & 0 \\
    0 & 0 \\
  \end{array}
\right]
  \begin{array}{l}
    \}r \\
    \}n-r \\
  \end{array}
.$$
Using the full matrix SVD of $A_{(j)}$,
$$A_{(j)}=U_{(j)}\Sigma_{(j)}V_{(j)}^*, \quad 2\leq j\leq d,$$
we have
\begin{align*}
(I_{n_j}-U_jU_j^*)A_{(j)} & = (U_{(j)}I_{n_j}U_{(j)}^*-U_{(j)}\tilde{I}_{r_j,n_j}U_{(j)}^*)U_{(j)}\Sigma_{(j)}V_{(j)}^* \\
& = U_{(j)}\left[
\begin{array}{cc}
  0 &  \\
   & I_{n_j-r_j} \\
\end{array}
\right]U_{(j)}^*U_{(j)}\Sigma_{(j)}V_{(j)}^* \\
& = U_{(j)}
\left[
  \begin{array}{cccccc}
    0 &  &  &  &  &  \\
     & \ddots &  &  &  &  \\
     &  & 0 &  &  &  \\
     &  &  & \sigma_{r_j+1}(A_{(j)}) &  &  \\
     &  &  &  & \ddots &  \\
     &  &  &  &  & \sigma_{n_j}(A_{(j)}) \\
  \end{array}
\right]
V_{(j)}^*.
\end{align*}
This implies that
\begin{equation}\label{tm:errorHOSVD}
\|(I_{n_j}-U_jU_j^*)A_{(j)}\|_F^2 = \sum_{i=r_j+1}^{n_j} \sigma_i^2(A_{(j)})
\leq (n_j-r_j)\sigma_{r_j+1}^2(A_{(j)}).
\end{equation}

Further on, using~\eqref{rel:CX} for $j=1$ and the property of the Moore-Penrose inverse,
$$CC^+C=C,$$
we get
\begin{align*}
(I_{n_1}-CC^+)A_{(1)} & = (I_{n_1}-CC^+)(CF+E) \\
& = CF+E-CC^+CF-CC^+E \\
& = (I_{n_1}-CC^+)E.
\end{align*}
Since $I_{n_1}-CC^+$ is projection we obtain
\begin{equation}\label{tm:PQR1}
\|(I_{n_1}-CC^+)A_{(1)}\|_F^2 \leq \|E\|_F^2=\|R_{22}\|_F^2,
\end{equation}
where $R_{22}$ is as in~\eqref{rel:QR} with $j=1$. Equality $\|E\|_F=\|R_{22}\|_F$ follows from~\eqref{rel:CXerror}.
To get the upper bound on the norm of $R_{22}$, we use~\cite[Lemma 2.5]{ACGri20} that gives
\begin{equation}\label{rel:R22acg}
\|R_{22}\|_2\leq\sqrt{1+2r_1+\sum_{j=1}^{r_1-1}4^j(r_1-j)}\sqrt{n_1-r_1}\sigma_{r_1+1}(A_{(1)})=\sqrt{p(r_1,n_1)}\sigma_{r_1+1}(A_{(1)}),
\end{equation}
where function $p$ is defined in~\eqref{def:p}.
Applying the equivalence of norms on~\eqref{rel:R22acg} it follows that
$$\|R_{22}\|_F\leq \sqrt{p(r_1,n_1)}\sqrt{n_1-r_1}\sigma_{r_1+1}(A_{(1)}),$$
that is
\begin{equation}\label{rel:R22}
\|R_{22}\|_F^2\leq p(r_1,n_1)(n_1-r_1)\sigma_{r_1+1}^2(A_{(1)}).
\end{equation}
We now use~\eqref{tm:PQR1} and~\eqref{rel:R22} to get
\begin{equation}\label{tm:errorPQR}
\|(I_{n_1}-CC^+)A_{(1)}\|_F \leq p(r_1,n_1)(n_1-r_1)\sigma_{r_1+1}^2(A_{(1)}).
\end{equation}

Finally, we insert relations~\eqref{tm:errorHOSVD} and~\eqref{tm:errorPQR} in~\eqref{tm:error1} to obtain the bound~\eqref{rel:tmerror}.
\end{proof}

Analogous bound as in Theorem~\ref{tm:error} can be obtained for the case~\eqref{def:hybrid-more} where we want to keep the original fibers in $t$ modes of a tensor. Then we get
\begin{equation}\label{def:hybrid-moreerror}
\|\calE\|_F^2\leq \sum_{i=1}^t p(r_i,n_i)(n_i-r_i)\sigma_{r_i+1}^2(A_{(i)}) + \sum_{j=t+1}^d (n_j-r_j)\sigma_{r_j+1}^2(A_{(j)}),
\end{equation}
with $p(r_i,n_i)$ defined in~\eqref{def:p}, $1\leq i\leq t$.

Contrary to the hybrid approach, assume that the approximation of $\calA$ is given as
$$\hat{\calA}_{\text{CUR}}=\calS\times_1C_1\times_2C_2\times_3\cdots\times_dC_d,$$
where matrices $C_i$, $1\leq i\leq d$, are obtained by PQR decomposition~\eqref{rel:C}.
Using the same reasoning it is easy to see that the error of such  approximation equals
\begin{equation}\label{rel:CURerror}
\|\calA-\hat{\calA}_{\text{CUR}}\|_F \leq \sum_{i=1}^d p(r_i,n_i)(n_i-r_i)\sigma_{r_i+1}^2(A_{(i)}),
\end{equation}
where $p(r_i,n_i)$ is as in~\eqref{def:p} for $1\leq i\leq d$.
Since $p(r,n)>1$, the difference between the error bounds~\eqref{rel:tmerror} and~\eqref{rel:CURerror} is increasing as the tensor order $d$ is increasing. Also, the difference between~\eqref{def:hybrid-moreerror} and~\eqref{rel:CURerror} is increasing as the number of modes $t$ in which we preserve the original fibers is decreasing.

\section{Numerical examples}\label{sec:numerical}

To illustrate the advantages of the hybrid CUR approximation we present three numerical examples. The tests are preformed using Matlab 2019b.

\begin{figure}[h]
    \centering
    \begin{subfigure}[b]{0.49\textwidth}
        \includegraphics[width=\textwidth]{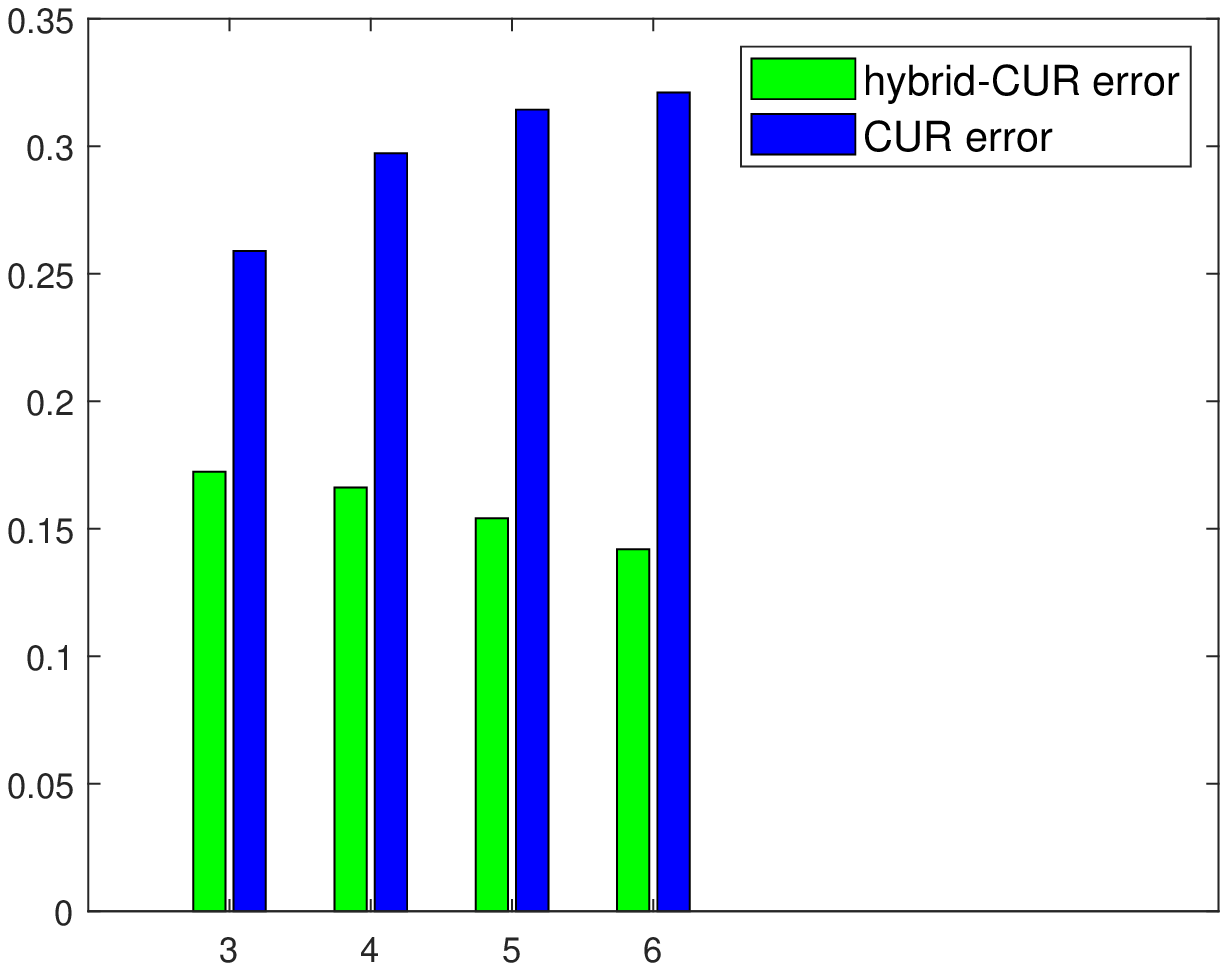}
        \caption{Tensor $\calA$}
    \end{subfigure}
    \begin{subfigure}[b]{0.49\textwidth}
        \includegraphics[width=\textwidth]{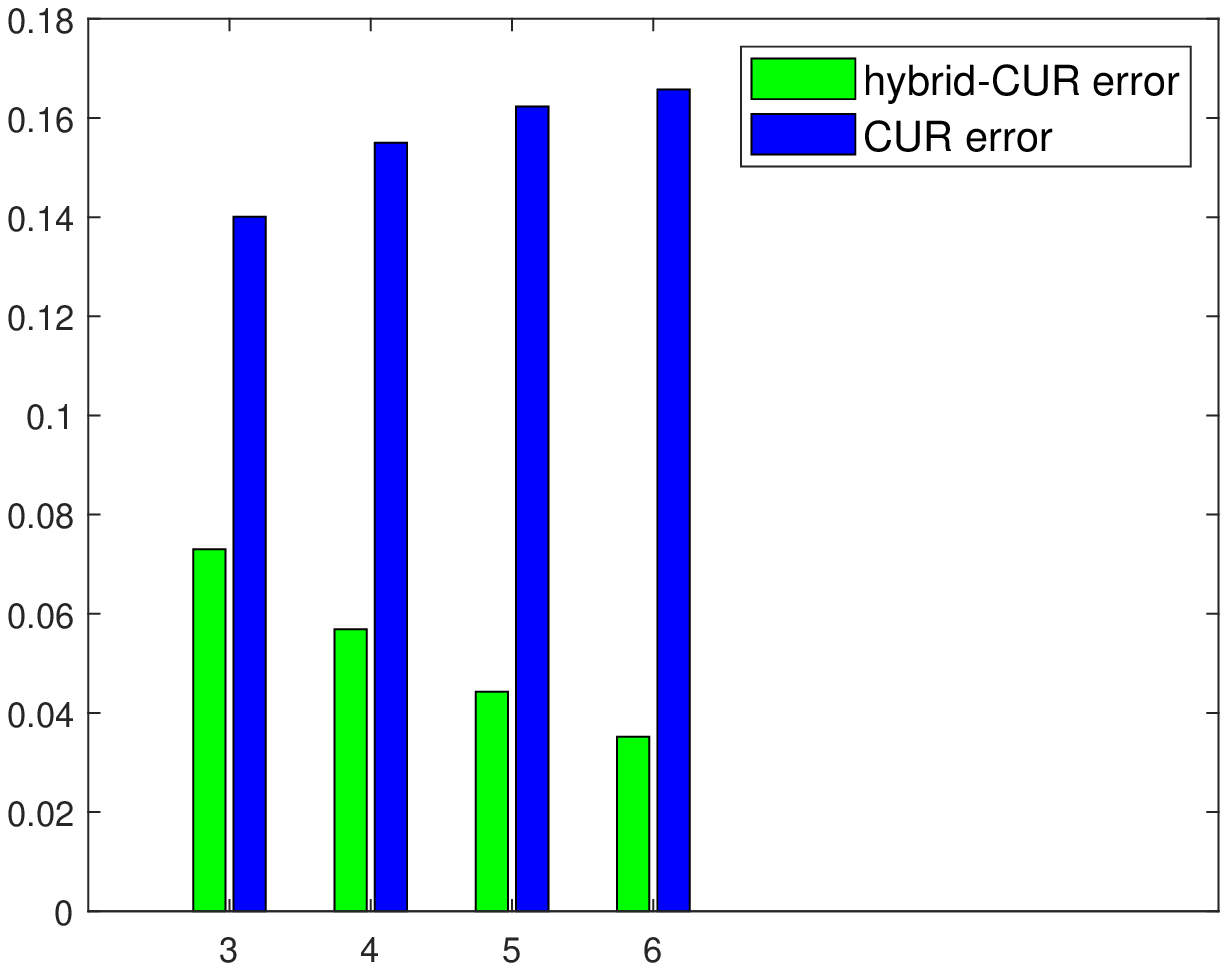}
        \caption{Tensor $\calB$}
    \end{subfigure}
    \caption{Relative approximation error when tensor order varies.}\label{fig:error-d}
\end{figure}

In Figure~\ref{fig:error-d} we compare the relative error obtained by HOID method from Algorithm~\ref{alg:hoid} and the hybrid method from Algorithm~\ref{alg:hybrid}. We show the results of the experiments performed on two function related tensors,
$$\calA(i_1,\ldots,i_d)=\frac{1}{i_1+i_2+\cdots+i_d}, \qquad \calB(i_1,\ldots,i_d)=\frac{1}{i_1+2\cdot i_2+\cdots+d\cdot i_d},$$
approximated with rank-$1$ tensors.
Here, tensor $\calA$ is symmetric, while $\calB$ is not. We set $n_1=n_2=\cdots=n_d=7$ and vary tensor order for $d=3,4,5,6$.
We observe that the relative error is significantly smaller when hybrid method is used. The difference between the relative errors is increasing with the tensor order $d$.

\begin{figure}[h]
    \centering
    \begin{subfigure}[b]{0.49\textwidth}
        \includegraphics[width=\textwidth]{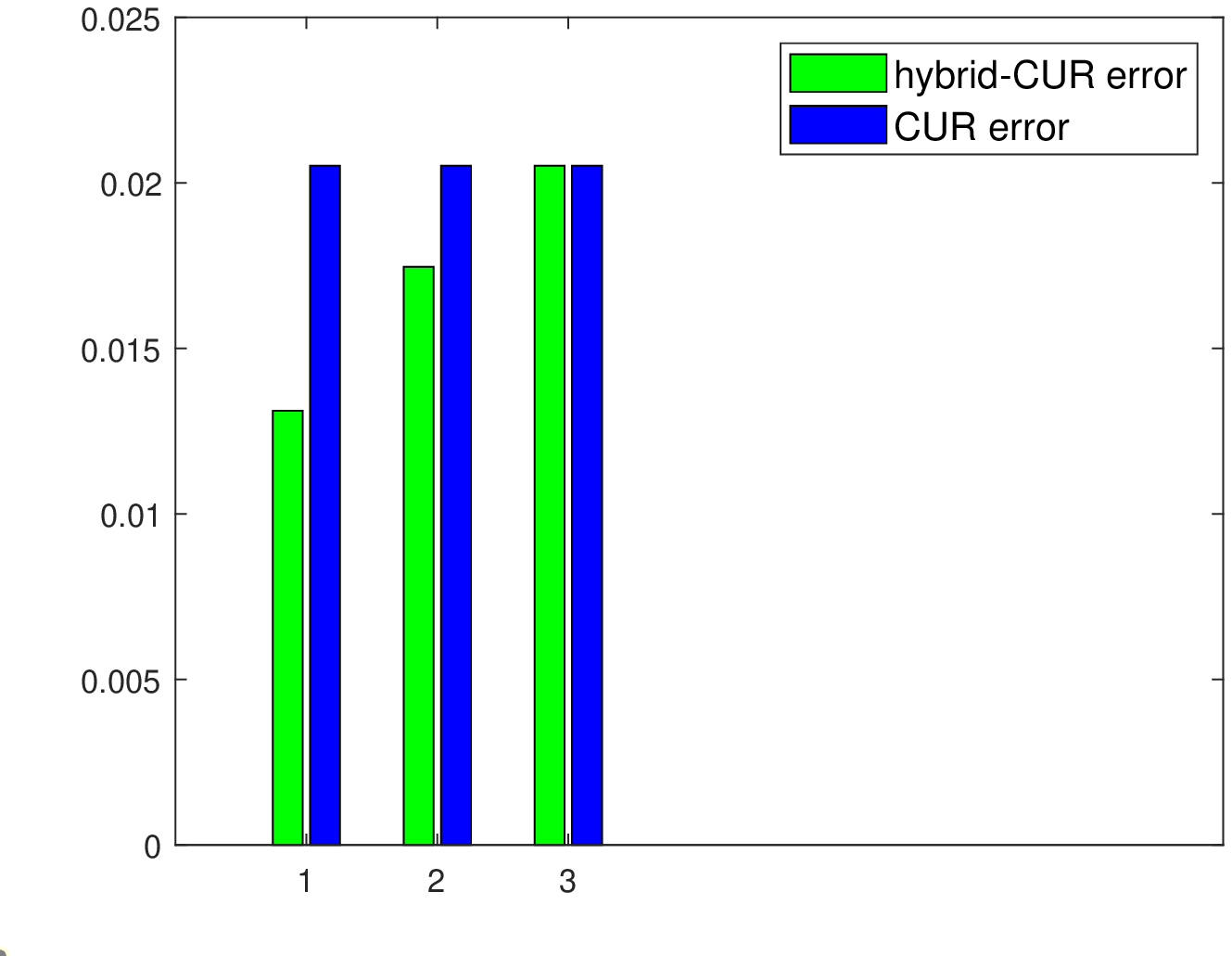}
        \caption{Tensor $\calB$ for $d=3$}
    \end{subfigure}
    \begin{subfigure}[b]{0.49\textwidth}
        \includegraphics[width=\textwidth]{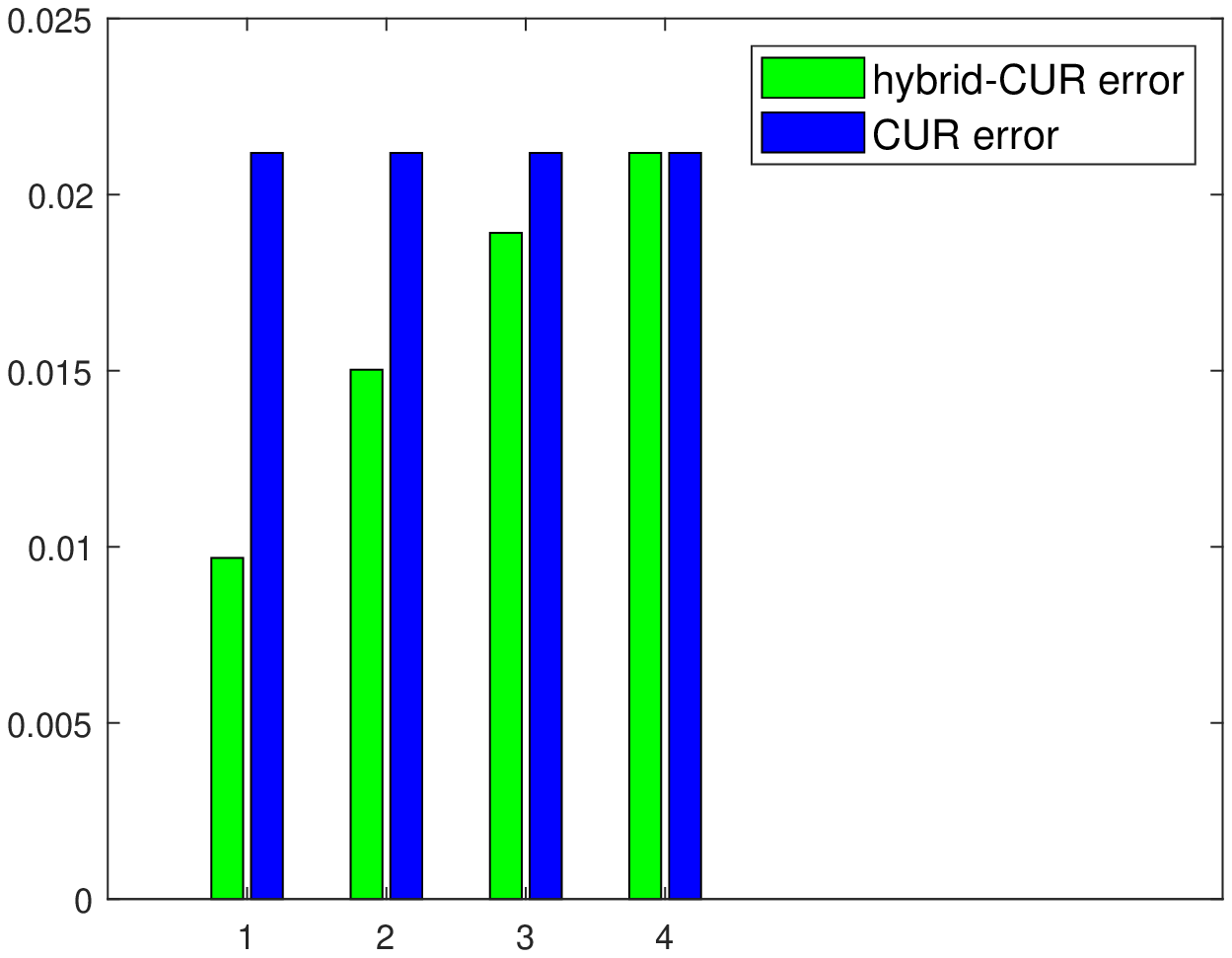}
        \caption{Tensor $\calB$ for $d=4$}
    \end{subfigure}
    \caption{Relative approximation error when the number of modes in which the fibers of the original tensors are preserved varies.}\label{fig:error-n}
\end{figure}

In Figure~\ref{fig:error-n} we also compare the relative error arising from Algorithm~\ref{alg:hoid} and Algorithm~\ref{alg:hybrid}. Here we vary the number of modes in which we preserve the original fibers when using the hybrid CUR. We perform the test on tensor $\calB$ for $d=3$ and $d=4$. We set $n_1=n_2=\cdots=n_d=30$ and do rank-$2$ approximation.
As expected, the difference in the approximation error is bigger as the number of the modes in which we preserve the original fibers is smaller. If we want to preserve the original fibers in all modes, than hybrid method boils down to the regular CUR method.

\begin{figure}[h]
    \centering
    \begin{subfigure}[b]{0.49\textwidth}
        \includegraphics[width=\textwidth]{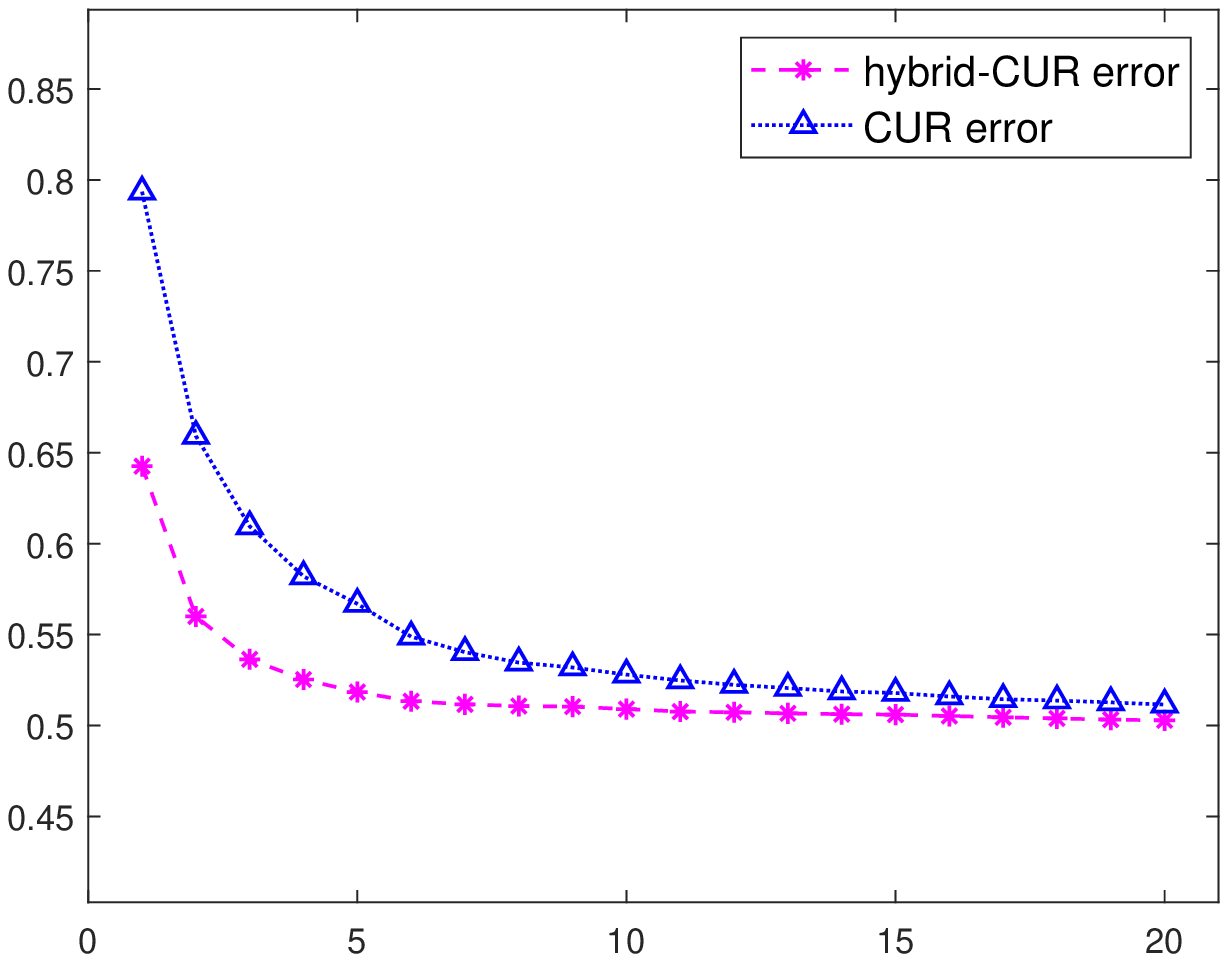}
        \caption{$d=3$, $n_1=n_2=n_3=100$}
    \end{subfigure}
    \begin{subfigure}[b]{0.49\textwidth}
        \includegraphics[width=\textwidth]{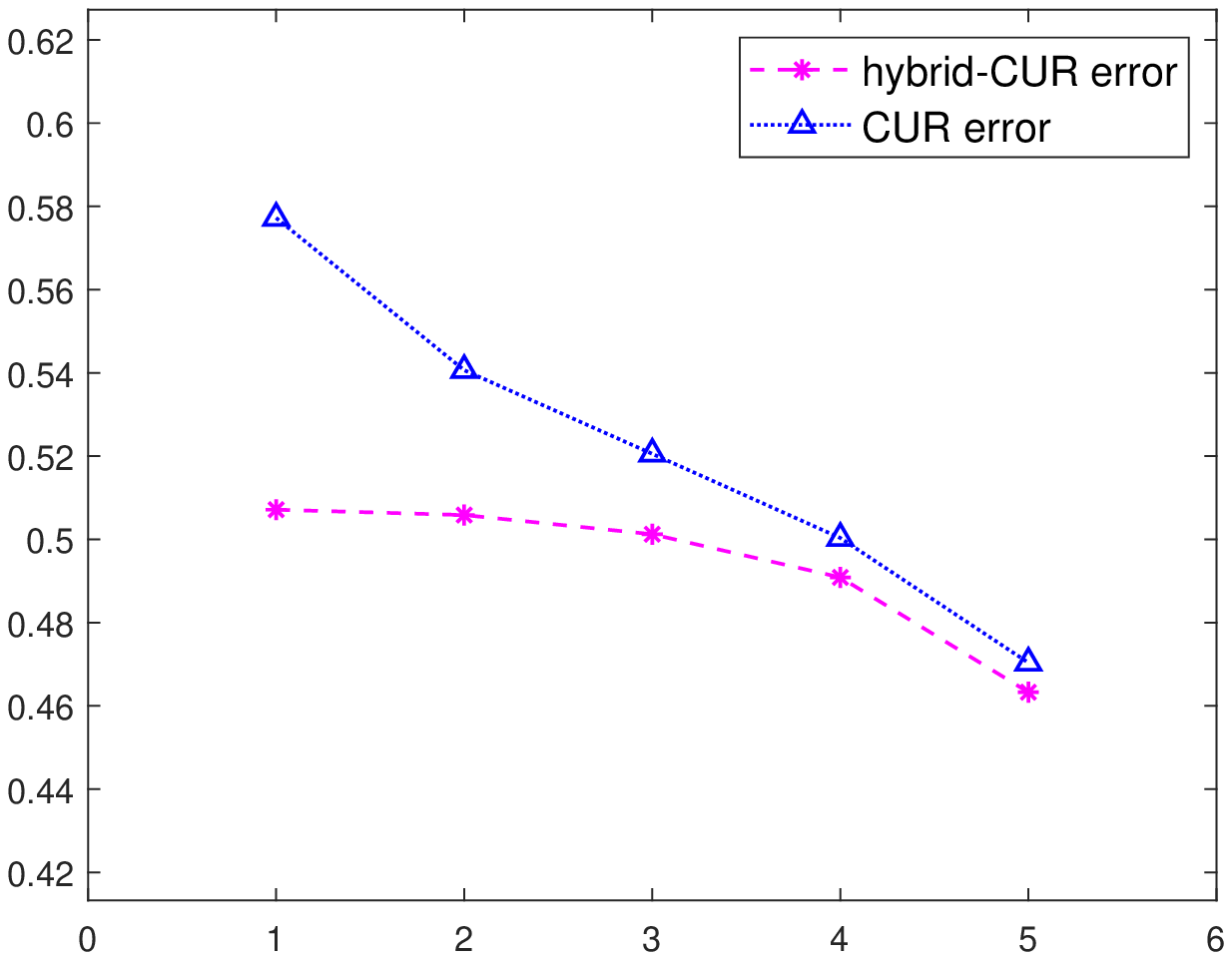}
        \caption{$d=6$, $n_1=\cdots=n_6=7$}
    \end{subfigure}
    \caption{Relative approximation error for random tensors when the approximation rank varies.}\label{fig:error-k}
\end{figure}

Moreover, in Figure~\ref{fig:error-k} we show the approximation error when the approximation rank $k$ varies. In the hybrid method we preserve the original fibers only in the first mode. We do the test for $d=3$ on a random $100\times100\times100$ tensor, and for $d=6$ on a random $7\times7\times7\times7\times7\times7$ tensor.

\section{Matrix case}\label{sec:matrix}

Hybrid CUR-type rank-$k$ approximation of a matrix $A\in\mathbb{R}^{m\times n}$ is a spacial case of~\eqref{def:hybrid-approx}.
Using tensor terminology, matrix columns are mode-$1$ fibers, while its rows are mode-$2$ fibers. For the ``matricizations'' of $A$ we have
$$A_{(1)}=A \quad \text{and} \quad A_{(2)}=A^T.$$

Assume that the approximation preserves the columns of $A$.
Then
\begin{equation}\label{def:hybrid-mcol}
\hat{A}=S\times_1C\times_2V=(CS)\times_2V=\left(V(CS)^T\right)^T=CSV^T,
\end{equation}
where $C\in\mathbb{R}^{m\times k}$ is made of $k$ columns of $A$. Matrix $V\in\mathbb{R}^{k\times n}$ consists of $k$ leading right singular vectors of $A$ attained by SVD, which are actually left singular vectors of $A_{(2)}$. Core ``tensor'' is matrix
$$S=A\times_1C^+\times_2V^T=C^+AV.$$

On the other hand, assume that the approximation preserves the rows of $A$. Rows of $A=A_{(1)}$ can also be considered as columns of $A_{(2)}$.
Here we have
$$\hat{A}=S\times_1U\times_2C=(US)\times_2C=\left(C(US)^T\right)^T=USC^T=USR,$$
where $R=C^T$ contains $k$ rows of $A$, $U$ contains $k$ leading left singular vectors of $A$, and
$$S=A\times_1U^T\times_2C^+=A\times_1U^T\times_2(R^T)^+=U^TAR^+.$$

Error obtained this way is smaller than the error obtained by CUR decomposition. This difference is quantified in Corollary~\ref{tm:error-matrix}. Its proof follows as a special case of Theorem~\ref{tm:error}.

\begin{Corollary}\label{tm:error-matrix}
Let $A\in\mathbb{C}^{m\times n}$. Let $\hat{A}$ be a rank-$k$ approximation of $A$ as in relation~\eqref{def:hybrid-mcol}. Assume that the matrix $C$ is obtained by QR with column pivoting. Then the approximation error $E$ satisfies the following inequality,
$$\|E\|_F=\|A-\hat{A}\|_F \leq p(k,m)(m-k)\sigma_{k+1}^2(A) + (n-k)\sigma_{k+1}^2(A),$$
where
$p(k,m)$ is defined by relation~\eqref{def:p}.
\end{Corollary}

\begin{figure}[h]
    \centering
        \includegraphics[width=0.5\textwidth]{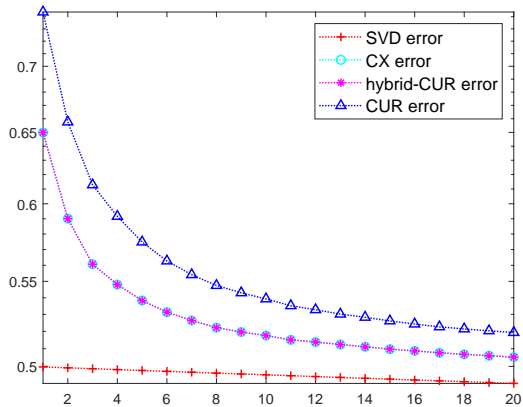}
    \caption{Relative approximation error for a random matrix when the approximation rank varies.}\label{fig:error-matrixcase}
\end{figure}

The rank-$k$ approximation error obtained by truncated SVD~\eqref{svdapproximation} is
$$\|E\|_2\geq\|A-A_k\|_2=\sigma_{k+1} \quad \text{and} \quad \|E\|_F\geq\|A-A_k\|_F=\sqrt{\sigma_{k+1}^2+\cdots+\sigma_n^2},$$
where $\sigma_i$ denotes the $i$-th singular value of $A$.
In Figure~\ref{fig:error-matrixcase} we illustrate the claim of Proposition~\ref{tm:error-matrix} by comparing the relative approximation error in the Frobenius norm obtained by rank-$k$ approximation of a matrix using four matrix decompositions: SVD (for the reference case), CX (where matrix $X$ is obtained as $X=C^+A$), hybrid CUR as in~\eqref{def:hybrid-mcol}, and matrix CUR. Test is done on a random $2000\times2000$ matrix. Approximation rank $1\leq k\leq20$ is given on the horizontal axis. In the matrix case, our hybrid CUR approach is equivalent to CX matrix decomposition.

\section*{Acknowledgements}
This work has been supported in part by Croatian Science Foundation under the project UIP-2019-04-5200.
The author would like to thank Georgia Tech for the kind hospitality during the process of writing this paper.

\end{document}